\documentclass[11pt]{amsart}
\usepackage{a4wide}
\usepackage{amsmath, amsfonts, amssymb,amscd,graphicx}
\usepackage{mathrsfs,url}

\date{\today}

\author
{Manuel Bodirsky}
\address{Institut f\"{u}r Algebra\\TU Dresden\\01062 Dresden\\Germany}
\email{Manuel.Bodirsky@tu-dresden.de}
   \urladdr{http://www.math.tu-dresden.de/~bodirsky/}
    \thanks{Both authors have received funding from the European Research Council under the European Community's Seventh Framework Programme (FP7/2007-2013 Grant Agreement no. 257039). The second author has been supported by funding of the Excellence Initiative by the German Federal and State Governments.}

\author
{Friedrich Martin Schneider}
    \address{Institut f\"{u}r Algebra\\TU Dresden\\01062 Dresden\\Germany}
    \email{Martin.Schneider@tu-dresden.de}
    \urladdr{http://tu-dresden.de/Members/friedrich\_martin.schneider/}
 
\title[Endomorphism Monoids of Countable Structures]{A topological characterisation of endomorphism monoids of countable structures}

\DeclareMathOperator{\ar}{ar}

\newcommand{\aC}{{\bf C}}

\newcommand{\sC}{\ensuremath{\mathscr{C}}}

\newcommand{\mN}{\mathbb N}

\newcommand{\ignore}[1]{}

\DeclareMathOperator{\Aut}{Aut}
\DeclareMathOperator{\Sym}{Sym}
\DeclareMathOperator{\End}{End}
\DeclareMathOperator{\Pol}{Pol}

\theoremstyle{plain}

    \newtheorem{thm}{Theorem}
    \newtheorem{theorem}[thm]{Theorem}

\theoremstyle{definition}

    \newtheorem{definition}{Definition}

\DeclareMathOperator{\comp}{comp}

\begin{document}
\maketitle
\begin{abstract}
A topological monoid is isomorphic to an
endomorphism monoid of a countable structure if and only if it is separable
and has a compatible complete ultrametric such that composition from the left 
is non-expansive. We also give a 
topological characterisation of those topological monoids that are isomorphic to
endomorphism monoids of countable $\omega$-categorical structures. 
Finally we present analogous characterisations 
for polymorphism clones of countable structures and for polymorphism clones of countable $\omega$-categorical structures. 
\end{abstract}

\section{Introduction}
\label{sect:intro}
Permutation groups and, more generally, transformation 
monoids on a set $X$ carry a natural topology such that composition is continuous, the \emph{topology of pointwise convergence}, which is induced by the product topology on $X^X$,
the set of all functions from $X$ to $X$, 
 where $X$ is considered to be discrete. 

A topological group is called 
\emph{non-archimedian} if it has a base at the identity consisting of open subgroups. 
The following is due to Becker and Kechris~\cite{BeckerKechris}.
\begin{thm}[Section 1.5 in~\cite{BeckerKechris}; also see Theorem~2.4.1 and Theorem 2.4.4 in~\cite{Gao}]\label{thm:groups}
\label{satz:topo-auto-groups}
Let $\bf G$ be a topological group. Then the following are equivalent.
\begin{enumerate}
\item $\bf G$ is topologically isomorphic to 
the automorphism group of a countable structure.  
\item $\bf G$ is topologically isomorphic to a closed subgroup of the full symmetric group $\Sym(\mathbb N)$. 
\item $\bf G$ is Polish and admits a compatible left invariant ultrametric.
\item $\bf G$ is Polish and non-archimedian. 
\end{enumerate}
\end{thm}
It has been asked in~\cite{Reconstruction} (Question 1 on page 29), and, independently, in~\cite{Tarzi} (`open question' on page 215), whether this theorem can 
be generalised from permutation groups to transformation monoids. Here, we present such a generalisation (Theorem~\ref{thm:monoids}):
a topological monoids $\bf M$ is topologically isomorphic to 
a closed submonoid of the full transformation monoid ${\mathbb N}^{\mathbb N}$ if and only if 
$\bf M$ is separable and admits a compatible complete left non-expansive ultrametric (Theorem~\ref{thm:monoids}). 

There is a strong link between closed subgroups of $\Sym(X)$ for countable $X$, and the model theory of structures $\Gamma$ with domain $X$.
This link is particularly strong 
when $\Gamma$ is \emph{$\omega$-categorical},
this is, when the first-order theory of $\Gamma$
has only one countable model up to isomorphism. 
These structures arise in many areas of mathematics; we refer to a recent survey article of MacPherson~\cite{MacphersonSurvey}. 
By the theorem of Ryll-Nardzewski, $\Gamma$ is $\omega$-categorical if and only if the automorphism group $\Aut(\Gamma)$ of $\Gamma$ is
\emph{oligomorphic}, this is, for all $n \in {\mathbb N}$ the action of $\Aut(\Gamma)$ on $X^n$ has 
only finitely many orbits. In fact, two 
$\omega$-categorical structures have 
topologically isomorphic automorphism groups
if and only if they are bi-interpretable, due to Coquand (see~\cite{AhlbrandtZiegler}). 
Hence, one might expect that the
property of a topological group to be topologically
isomorphic to a closed oligomorphic permutation group can be expressed directly in terms of the
topological group (without reference to any action). 
This is indeed possible, and closely related
to a result of Todor Tsankov about 
closed subgroups of $\Sym({\mathbb N})$ that are Roelcke precompact~\cite{Tsankov}. We will review the characterisation of topological groups with an oligomorphic action on a countable set in Section~\ref{sect:ogroups}. 

When two $\omega$-categorical structures
$\Gamma$ and $\Delta$ 
have not only topologically isomorphic automorphism groups, but even topologically isomorphic endomorphism monoids, one might expect a stronger form of reconstruction than
reconstruction up to bi-interpretability.
And indeed, if $\Gamma$ and $\Delta$
do not have constant endomorphisms, then
$\End(\Gamma)$ and $\End(\Delta)$
are topologically isomorphic if and only if 
$\Gamma$ and $\Delta$ are existentially bi-interpretable~\cite{BodJunker}.
In Section~\ref{sect:oligo-monoids} we present
a characterisation of the topological monoids
that appear as endomorphism monoids of countable 
$\omega$-categorical structures, and this specialises
to the mentioned result for oligomorphic permutation groups. 

We finally generalise our results further to function clones (Theorem~\ref{thm:clones}).
A subset of ${\mathscr O}_X := \bigcup_{n \in \mathbb N} (X^n \to X)$
is called a \emph{function clone} if it contains the projections and is closed under composition. 
\emph{Clones} are abstractions of function clones,
analogously as groups are abstractions of permutation groups, 
and monoids are abstractions of transformation monoids. It is well known that every abstract clone can be realised as a function clone. 

For a set $X$ and $k \in \mathbb N$, we equip $X^k \to X$ with the product topology where $X$ is taken to be discrete,
and we take $\bigcup_{k \in \mathbb N} X^k \to X$ as the sum space. Observe that with respect to this topology, composition in function clones is continuous.
In Section~\ref{sect:clones},
we characterise the 
topological clones that are isomorphic to closed subclones of ${\mathscr O}_{X}$. 
These are precisely the clones that arise as \emph{polymorphism clones} 
$\Pol(\Gamma)$ 
of countable structures $\Gamma$.
\emph{Polymorphisms} are multivariate endomorphisms and an important
concept from universal algebra that lately has
seen many applications in 
theoretical computer science (see, e.g., the collection of survey articles~\cite{CSPSurveys}).

When a countable structure $\Gamma$ is
$\omega$-categorical, then the polymorphism clone $\Pol(\Gamma)$, viewed as a topological clone, determines $\Gamma$ up
to \emph{primitive positive} bi-interpretability~\cite{Topo-Birk}. Here, in contrast to the result on monoids mentioned above, one does not need the assumption 
that $\Pol(\Gamma)$ and $\Pol(\Delta)$ do not have constant operations. 
Primitive positive interpretability plays a central role for
the study of the complexity of constraint satisfaction problems. In particular, this result
implies that the computational complexity of 
the constraint satisfaction problem of $\Gamma$
only depends on $\Pol(\Gamma)$, viewed
as a topological clone. In Section~\ref{sect:oligo-clones}, we characterise the topological clones
that can arise as polymorphism clones of $\omega$-categorical structures directly in terms of the topology. 

\begin{figure}
\begin{center}
\begin{tabular}{|l||l|l|}
\hline
& General case & Oligomorphic case \\
\hline
\hline
Groups & Theorem~\ref{thm:groups} (Becker-Kechris) & Theorem~\ref{thm:ogroups} (essentially Tsankov) \\
& Section~\ref{sect:intro} & Section~\ref{sect:ogroups} \\
\hline
Monoids & Theorem~\ref{thm:monoids} (Question of Tarzi) & Theorem~\ref{thm:omonoids} \\
& Section~\ref{sect:monoids} & Section~\ref{sect:oligo-monoids} \\
\hline
Clones & Theorem~\ref{thm:clones} & Theorem~\ref{thm:oclones} (capture CSP complexity) \\
& Section~\ref{sect:clones} & Section~\ref{sect:oligo-clones} \\
\hline
\end{tabular}
\end{center}
\caption{Overview of the results.}
\end{figure}

\section*{Notation}
We denote the set of natural numbers by ${\mathbb N} = \{ 1,2,\ldots \}$. Given any $n \in {\mathbb N}$, we write $\underline{n}$ for the set $\{1,\dots,n\}$. If $f \colon X \to Y$,
$t \in X^k$, then we write $f(t)$ for the
tuple $(f(t_1),\dots,f(t_k)) \in Y^k$. 
For an equivalence relation $U$ on $X$, and $x \in X$, we write
$[x]_U$ for the equivalence class of $x$ with respect to $U$, and $X/U$ for the set $\{[x]_U  \; | \; x \in X\}$.
The cardinality of $X/U$ is called the \emph{index} of $U$.

\section{Monoids}
\label{sect:monoids}
Let ${\bf M} = (M;{\cdot},1)$ be a topological monoid and $d$ a compatible
metric on $M$. 
We say that $d$ 
is \emph{left non-expansive} if for all $g,f,f' \in M$ we have that
$d(gf,gf') \leq d(f,f')$. A metric $d$ on a set $X$ 
is called an \emph{ultrametric} if $d(x,z) \leq \max(d(x,y),d(y,z))$ for all $x,y,z \in X$. 
In this section we show that 
$\bf M$ is \emph{topologically isomorphic} to a closed submonoid of $\mN^{\mN}$ (that is, there is a monoid isomorphism which is a homeomorphism)  if and only if $\bf M$ is separable and admits a compatible complete left non-expansive ultrametric. 

In the proof we need the following fundamental concepts from topology (see e.g.~\cite{Bourbaki}). 
A nonempty collection $\mathcal U$ 
of binary reflexive relations on a set $X$ 
is a \emph{uniformity} (and $X$ along with $\mathcal U$ is called a
\emph{uniform space}) 
if it satisfies the following axioms:
\begin{itemize}
\item If $U \in \mathcal U$ and $U\subseteq V \subseteq X^2$, then $V \in \mathcal U$.
\item If $U \in \mathcal U$, then $U^{-1} \in \mathcal U$.
\item If $U \in \mathcal U$ and $V \in \mathcal U$, then $U\cap V \in \mathcal U$.
\item If $U \in \mathcal U$, then there is $V \in \mathcal U$ such that $V \circ V \subseteq U$. 
\end{itemize}
The elements of $\mathcal U$ are called \emph{entourages}. 
A \emph{base} of a uniformity $\mathcal{U}$ is a subset ${\mathcal B} \subseteq {\mathcal U}$ 
such that
for any $V \in \mathcal U$ there is a $W \in {\mathcal B}$ such that $W \subseteq V$. Clearly, any base $\mathcal B$ 
of a uniformity $\mathcal U$ determines 
$\mathcal U$ uniquely, and we call $\mathcal U$ 
the uniformity \emph{generated} by $\mathcal B$. 
The topology \emph{induced} by $\mathcal U$ on $X$
is the topology where $O \subseteq X$ is open if for every $x \in O$ there exists an entourage $U \in \mathcal U$
such that $\{y \; | \; (x,y) \in U\} \subseteq O$. 
In this case we say that $\mathcal U$ is \emph{compatible}
with the topology. 
Every metric space $(X,d)$ \emph{generates} 
a uniformity, namely the uniformity with the base 
$\mathcal U = \big \{ \{(x,y) \in X^2 : d(x,y) \leq \epsilon\} \mid \epsilon \in {\mathbb R} \big \}$.  


A \emph{Cauchy filter} on a uniform space $X$
is a filter $\mathcal F$ on $X$
such that for every $U \in \mathcal U$ there exists an $A \in \mathcal F$ with $A \times A \subseteq U$. We say that $\mathcal F$ \emph{converges} if there exists an $x \in X$ such that for every $U \in \mathcal U$ we 
have that $\{y \; | \; (x,y) \in U\} \in \mathcal F$, in which case we say that \emph{$\mathcal F$ converges to $x$}. 
The uniform space $X$ is called \emph{complete} if every Cauchy filter on $X$ converges.
A function $f \colon X \to Y$ between uniform spaces is called \emph{uniformly continuous} if the pre-image of
any entourage in $Y$ is an entourage in $X$. 

\begin{thm}\label{thm:monoids}
\label{satz:topo-transformation-monoids}
Let $\bf M$ be a topological monoid. Then the following are equivalent.
\begin{enumerate}
\item\label{m-struct} $\bf M$ is topologically isomorphic to 
the endomorphism monoid of a countable structure.  
\item\label{m-start} $\bf M$ is topologically isomorphic to a closed submonoid of $\mN^{\mN}$. 
\item\label{m-goal} $\bf M$ is separable and admits a compatible complete left non-expansive ultrametric.
\item\label{m-uniformity} $\bf M$ is Hausdorff and has a 
complete compatible uniformity with a countable base $\mathcal B$ such that 
\begin{enumerate}
\item\label{a} each element of $\mathcal B$ is an equivalence relation on $M$ of countable index, and
\item\label{b} for all $U \in \mathcal B$ and for all $x,y,z \in M$, if $(y,z) \in U$ then $(xy,xz) \in U$. 
\end{enumerate}
\end{enumerate}
\end{thm}

\begin{proof}
The equivalence of $(\ref{m-struct})$ and $(\ref{m-start})$ is well-known
and can 
e.g.\ be found in~\cite{Bodirsky-HDR} (Proposition 3.4.8). We show the implications $(\ref{m-start}) \Rightarrow (\ref{m-goal})$, 
$(\ref{m-goal}) \Rightarrow (\ref{m-uniformity})$, 
$(\ref{m-uniformity}) \Rightarrow (\ref{m-start})$.

$(\ref{m-start}) \Rightarrow (\ref{m-goal})$.
Closed submonoids $\bf M$ 
of $\mN^{\mN}$ are separable: for all $s,t \in \mN^n$ such that there
exists an $f \in M$ with $f(s) = t$, we arbitrarily select such an $f$; the set of all selected $f$ lies dense in $M$. A compatible complete non-expansive ultrametric $d$ can be obtained as follows. 
For elements $f,g \in M$ we define $d(f,g) := 0$ if $f=g$,
and otherwise
$d(f,g) := 1/2^n$ where $n \in \mathbb N$ is least 
such that $f(n) \neq g(n)$. 

$(\ref{m-goal}) \Rightarrow (\ref{m-uniformity})$.
Let $d$ be a compatible left non-expansive ultrametric of $\bf M$. 
Then the collection ${\mathcal B} := \{B_n \; | \; n \in \mN\}$ where $B_n := \{(x,y) \in M^2  \; | \; d(x,y) \leq 1/n\}$ is a uniformity base which we claim
has the required properties. 
Obviously, ${\mathcal B}$ is countable,
and the transitivity of  the relations $B_n$ follows from the assumption that $d$ is an ultrametric. 
The uniformity generated by 
$\mathcal B$ is complete, because
it is generated by $d$, and 
the metric $d$ is complete. 
By separability of $M$, there exists 
a countable dense subset 
$\{e_1,e_2,\dots\}$ of $M$,
that is, for every $e \in M$ and $n \in \mathbb N$ 
 there exists an $i \in \mathbb N$ such that 
 $d(e_i,e) < 1/n$. This shows that $B_n$ has countably many equivalence classes. 
Finally, 
(\ref{b}) holds since the metric $d$ is left non-expansive. 

$(\ref{m-uniformity}) \Rightarrow (\ref{m-start})$. 
Let $V_1,V_2,\dots$ be an enumeration of the equivalence relations in $\mathcal B$, and 
define $(U_n)_{n \in \mathbb N}$ to be the sequence of equivalence relations given by 
$U_n := \bigcap_{i \in \underline n} V_i$.
Let $N_n := M / {U_n}$. Define $\xi \colon M \rightarrow N^N$ where $N$ is the countable set $\bigcup_{n \in \mathbb N} (N_n \times \{n\})$ as follows.
$$ \xi(x)([y]_{U_n},n) := ([xy]_{U_n},n)$$
This is well-defined: if $(y,y') \in U_n$, 
then $(xy,xy') \in U_n$ by~(\ref{b}). 
Moreover, $\xi$ is a homomorphism since
\begin{align*}
\xi(xy)([z]_{U_n},n) = ([xyz]_{U_n},n) = \xi(x)([yz]_{U_n},n) = 
\xi(x)(\xi(y)([z]_{U_n},n))
\end{align*}
and since $\xi(1)([z]_{U_n},n) = ([1z]_{U_n},n) = ([z]_{U_n},n)$.


\vspace{.2cm}
{\bf Claim 1.} The map $\xi$ is injective. 
For distinct $x,y \in M$ there exists an $n \in \mN$ such that $(x,y) \notin U_n$ since $M$ is Hausdorff.
Hence, $[x]_{U_n} \neq [y]_{U_n}$ and $[x1]_{U_n} \neq [y1]_{U_n}$. Therefore, $\xi(x)([1]_{U_n},n) \neq \xi(y)([1]_{U_n},n)$
and $\xi(x) \neq \xi(y)$. 


\vspace{.2cm}
{\bf Claim 2.} The map $\xi$ is continuous.
Let $x \in M$ be arbitrary, and $E \subseteq N$ be finite.
Define $W_{x,E} := \{g \in \xi(M)  \; | \; g(e) = \xi(x)(e) \text{ for all $e \in E$} \}$, and observe that the sets of the form $W_{x,E}$ form a 
basis for the topology induced by $\xi(M)$ in $N^N$.
Then 
\begin{align*}
\xi^{-1}(W_{x,E}) & = \{ y \in M \; | \; \xi(y)(e) = \xi(x)(e)\text{ for all $e \in E$} \} \\
& = \bigcap_{([z]_{U_n},n) \in E} \{ y \in M  \; | \; \xi(x)([z]_{U_n},n) = \xi(y)([z]_{U_n},n)\}  \\
& = \bigcap_{([z]_{U_n},n) \in E} \{ y \in M  \; | \; ([xz]_{U_n},n) = ([yz]_{U_n},n)\} \\
& = \bigcap_{([z]_{U_n},n) \in E} \{ y \in M  \; | \; yz \in [xz]_{U_n}\} \, .
\end{align*}
We will prove that $S:=\{ y \in M  \; | \; yz \in [xz]_{U_n}\}$ is open.  
Define $\rho_z \colon M \to M$
by $\rho_z(y) := y z$, which is continuous as the multiplication in $\bf M$ is continuous. Since $[xz]_{U_n}$ is open, 
it follows that $\rho_z^{-1}([xz]_{U_n}) = \{ y \in M \; | \; \rho_z(y) \in [xz]_{U_n}\} = S$
is open, too. Therefore, $\xi^{-1}(W_{x,E})$ is a finite intersection of open sets, and open. 
It follows that $\xi$ is continuous at $x$. 

\vspace{.2cm}
{\bf Claim 3.} The map $\xi^{-1}$ is uniformly continuous. We have to show that for all $n \in \mathbb N$ the relation
$\{(\xi(x),\xi(y))  \; | \; (x,y) \in U_n\}$ 
is  an entourage in $\xi(M)$. 
Define $$V_n := \{(f,g) \in N^N \times N^N  \; | \; f(([1_M]_{U_n},n)) = g(([1_M]_{U_n},n))\} \; .$$
Then 
\begin{align*}
V_n \cap \xi(M)^2 = & \; \big \{(\xi(x),\xi(y)) \; | \; x,y \in M \text{ and } \xi(x)([1_M]_{U_n},n) = \xi(y)([1_M]_{U_n},n) \big \} \\
= & \; \big \{(\xi(x),\xi(y)) \; | \; x,y \in M \text{ and } [x1_M]_{U_n} = [y1_M]_{U_n} \big \}  \\ 
= & \; \big \{(\xi(x),\xi(y)) \; | \; (x,y) \in U_n \big \} \, . 
\end{align*}
The former set is an entourage in $\xi(M)$, and this proves the claim. 

\vspace{.2cm}
{\bf Claim 4.} The image of $M$ under $\xi$ is closed in $N^N$. Let $(x_i)_{i \in \mN}$ be a sequence of
elements of $M$ such that $(\xi(x_i))_{i \in \mN}$ 
converges to $y \in N^N$. Let $d$ be
the metric on $N^N$ defined analogously to the metric on ${\mathbb N}^{\mathbb N}$ from the implication $(1) \Rightarrow (2)$. Then  
the set of all sets of the form $\{z \in N^N  \; | \; d(z,y) \leq 1/n\}$, for $n \in {\mathbb N}$, is a Cauchy filter $\mathcal F$.
 Since $\xi^{-1}$ is uniformly continuous, 
${\mathcal F}' := \{\xi^{-1}(A)  \; | \; A \in \mathcal F\}$ is a Cauchy filter, too. 
Since the uniformity generated by $\mathcal B$
is complete, ${\mathcal F}'$ converges 
to an element $x \in M$. 
Then $\xi(x)=y$ by continuity of $\xi$.
\end{proof}

\section{Oligomorphic Transformation Monoids}
\label{sect:oligo-monoids}
In this section we characterise those topological monoids
that appear as endomorphism monoids of $\omega$-categorical
structures. To this end, let us first recall the concept of \emph{approximate oligomorphicity} from \cite{Rosendal13} (see also \cite{YaacovTsankov}).

\begin{definition}\label{definition:approximately.oligomorphic}
	An isometric action of a group $G$ on a metric space $(X,d)$ is called \emph{approximately oligomorphic} if for any $n \geq 1$ and $\varepsilon > 0$ there is a finite subset $F \subseteq X^{n}$ such that \begin{displaymath}
		\forall x \in X^{n} \, \exists y \in F \colon \, d_{G,n}(x,y) < \varepsilon ,
\end{displaymath} where $d_{G,n}(x,y) := \inf_{g \in G} \max_{i \in \underline n}  d(x_i,g y_i)$ for $x,y \in X^{n}$. \end{definition}

Let $\bf M$ be a topological monoid. Consider the group $G := \{ x \in M \mid \exists y \in M \colon \, xy = yx = 1 \}$ of its \emph{units}. Note that if $d$ is a left non-expansive compatible metric on $M$, then the left action $G \times M \to M, \, (g,x) \mapsto gx$ is an isometric action of $G$ on $(M,d)$. Furthermore, let us introduce the following terminology.

\begin{definition}\label{def:g-finitely}
A metric $d$ on $M$ is called \emph{$G$-finitely generated}
if for every $\varepsilon > 0$ there 
exists a finite subset $F$ of $G$ such that
$$ \forall x,y \in M: \max_{g \in F} d(xg,yg) < 1/2
\Rightarrow d(x,y) < \varepsilon \; .$$
\end{definition}




\begin{thm}\label{thm:omonoids}
Let $\bf M$ be a topological monoid and let $G$
be the units of $\bf M$. Then the following are equivalent.
\begin{enumerate}
\item\label{om-struct} $\bf M$ is topologically isomorphic to 
the transformation monoid of a countable $\omega$-categorical structure.   
\item\label{om-start} 
There is a topological isomorphism $\xi$ between 
$\bf M$ and 
a closed submonoid of $\mN^{\mN}$ 
such that $\xi(G)$ is oligomorphic. 
\item\label{om-goal} $\bf M$ is separable and there exists a compatible $G$-finitely generated complete left non-expansive ultrametric $d$ on $M$ such that the left action of $G$ on $(M,d)$ is approximately oligomorphic.
\item\label{om-uniformity} $\bf M$ is Hausdorff and there exists a left invariant equivalence relation $V$ of countable index on $M$ such that
\begin{enumerate}
\item\label{om-uniformity-2} 
the uniformity $\mathcal U$ generated by $$\big \{ \{ (x,y) \in M^2 \mid \forall g \in F \colon \, (xg,yg) \in V \} \mid F \subseteq G \textit{ finite} \big \}$$ is compatible and complete, and
\item\label{om-uniformity-1} for each $n \in \mathbb{N}$, the equivalence relation \begin{displaymath}
		V_n := \{ (x,y) \in M^{n} \times M^{n} \mid \exists g \in G \, \forall i \in \underline{n} \colon \, (x_{i},gy_{i}) \in V \}
		\end{displaymath} has finite index.
\end{enumerate}
\end{enumerate}
\end{thm}

\begin{proof}
The equivalence of $(\ref{om-struct})$ and $(\ref{om-start})$ follows from the equivalence of $(\ref{om-struct})$ and $(\ref{om-start})$ in Theorem~\ref{thm:monoids}, and the theorem of Engeler, Svenonius, and Ryll-Nardzewski (see, e.g., Hodges~\cite{HodgesLong}). 
We show the implications $(\ref{om-start}) \Rightarrow (\ref{om-goal})$, 
$(\ref{om-goal}) \Rightarrow (\ref{om-uniformity})$, 
$(\ref{om-uniformity}) \Rightarrow (\ref{om-start})$. 

$(\ref{om-start}) \Rightarrow (\ref{om-goal})$.
It follows from the implication 
$(\ref{m-start}) \Rightarrow (\ref{m-goal})$ in Theorem~\ref{thm:monoids} that 
$\bf M$ is separable and admits a compatible complete left non-expansive ultrametric $\tilde d$.  
Since $\xi(G)$ is an oligomorphic permutation group, 
there exists some $n \in \mN$ such that $\mN = \bigcup  _{j \in \underline{n}} \{ \xi(g)(j) \mid g \in G\}$.
Then $d$ defined by
$d(x,y) := 2^n \tilde d(x,y)$ is also a compatible left non-expansive ultrametric on $M$. 
We first verify that $d$ is $G$-finitely generated. 
Let $\varepsilon > 0$. There exists $m \in \mN$ such that $2^{n-m} < \varepsilon$. For each $i \in \underline{m}$, choose some $g_{i} \in G$ such that $i \in \xi(g_i)(\underline{n})$. Evidently, $F := \{ g_{i} \mid i \in \underline{m} \}$ is finite. Let $x,y \in M$. 
If $\max_{g \in F} d(xg,yg) < 1/2$, then
$(\xi(xg_{i})){\restriction_{\underline{n}}} = (\xi(yg_{i})){\restriction_{\underline{n}}}$ for each $i \in \underline{m}$ and thus $\xi(x){\restriction}_{\underline{m}} = \xi(y){\restriction}_{\underline{m}}$, which implies that $\tilde d(x,y) \leq 2^{-m}$, and
hence $d(x,y) \leq 2^{n-m} < \varepsilon$. 

To complete the proof of $(\ref{om-goal})$, 
let $\varepsilon > 0$. 
Again choose $m \in \mN$ such that $2^{n-m} < \varepsilon$. Since $\xi(G)$ is oligomorphic, 
$\mathcal{Q} := \mN^{mn}/\xi(G)$ is finite. 
In the following, we identify the elements
of $M$ with their image under $\xi$ when
applying them to elements from $\mathbb N$. 
Let \begin{displaymath}
	\mathcal{Q}' := \{ Q \in \mathcal{Q} \mid \exists f \in M^n \colon \, (f_1(1),\ldots,f_1(m), \ldots \ldots ,f_n(1),\ldots,f_n(m)) \in Q \} .
\end{displaymath} For each $Q \in \mathcal{Q}'$, choose $f_{Q} \in M^{n}$ such that \begin{displaymath}
	(f_{Q,1}(1),\ldots,f_{Q,1}(m), \ldots \ldots ,f_{Q,n}(1),\ldots,f_{Q,n}(m)) \in Q .
\end{displaymath} Clearly, $F := \{ f_{Q} \mid Q \in \mathcal{Q}' \}$ is a finite subset of $M^{n}$. Let $f \in M^{n}$. Then there exists some $Q \in \mathcal{Q}'$ such that $(f_1(1),\ldots,f_{1}(m), \ldots \ldots ,f_{n}(1),\ldots,f_{n}(m)) \in Q$. 
By the choice of $f_Q$, 
there exists $g \in G$ such that \begin{align*}
	(f_1(1),\ldots, &f_1(m), \ldots \ldots ,f_{n}(1),\ldots,f_{n}(m)) \\
		&= (g(f_{Q,1}(1)),\ldots,g(f_{Q,1}(m)), \ldots \ldots ,g(f_{Q,n}(1)),\ldots,g(f_{Q,n}(m))) .
\end{align*} Consequently, for every $i \in \underline{n}$ 
we have $\tilde d(f_i,gf_{Q,i}) < 2^{-m}$ and hence 
$d(f_i,gf_{Q,i}) \leq 2^{n-m} < \varepsilon$. Therefore, $d_{G,n}((f_1,\dots,f_n),(f_{Q,1},\dots,f_{Q,n})) < \varepsilon$ and we are done.

$(\ref{om-goal}) \Rightarrow (\ref{om-uniformity})$.
Let $d$ be as in (\ref{om-goal}). 
Since $d$ is an ultrametric, $V := \{ (x,y) \in M^2 \mid d(x,y) < 1/2 \}$ is an equivalence relation on $M$. Besides, as $d$ is left non-expansive, $V$ is left invariant. Moreover, $M/V$ is countable because $M$ is separable. 
We now show (\ref{om-uniformity-2}).
By definition, the uniformity $\mathcal U$ 
is generated by 
$$\big \{ \{ (x,y) \in M^2 \mid \forall g \in F \colon \, d(xg,yg) < 1/2 \} \mid F \subseteq G \text{ finite} \big \} \; .$$
Since $d$ is $G$-finitely generated
we find for every $n \in \mathbb N$
a finite subset $F \subseteq G$ such that
$$ \forall x,y \in M : \max_{g \in F} d(xg,yg) < 1/2 \Rightarrow d(x,y) < 1/n \; .$$
It follows that $\mathcal U$ is generated
by 
$$\big \{ \{ (x,y) \in M^2 \mid d(x,y) < 1/n\} \mid n \in \mathbb N \big \}$$
and hence coincides 
with the uniformity generated by $d$. In particular, $\mathcal{U}$ is complete and compatible with the topology of $M$.

Furthermore, (\ref{om-goal}) implies (\ref{om-uniformity-1}). To see this, let 
$n \in \mN$. By (\ref{om-goal})
there exists a finite set $F \subseteq M^n$ such that for all $x \in M^n$ there exists a $y \in F$ 
such that 
$d_{G,n}(x,y) < 1/n$. 
Hence, the index of $V_n = \{(x,y) \in M^n \times M^n \; | \; 
d_{G,n}(x,y) < 1/n\}$ is bounded by $|F|$. 

\vspace{.2cm}
$(\ref{om-uniformity}) \Rightarrow (\ref{om-start})$.
By assumption, $N := M/V$ is countable. Let us define $\xi \colon M \rightarrow N^N$ by $$\xi(x)([y]_{V}) := [xy]_{V}$$ for all $x,y \in M$. Since $V$ is left invariant, $\xi$ is well defined. Furthermore, $\xi$ is a homomorphism because 
$$\xi(xy)([z]_{V}) = ([xyz]_{V}) = \xi(x)([yz]_{V}) = \xi(x)(\xi(y)([z]_{V}))$$ and 
$$\xi(1)([z]_{V}) = [1z]_{V} = [z]_{V}$$ for all $x,y,z \in M$.



\vspace{.2cm}
{\bf Claim 1.} The map $\xi$ is injective. 
Let $x,y \in M$ where $x \ne y$. As a consequence of (\ref{om-uniformity-2}) and $M$ being Hausdorff, there exists some $g \in G$ such that $(xg,yg) \notin V$.
Accordingly, $\xi(x)([g]_{V}) = [xg]_{V} \ne \xi(y)([g]_{V}) = [yg]_{V}$
and therefore $\xi(x) \neq \xi(y)$. 

\vspace{.2cm}
{\bf Claim 2.} $\xi(G)$ is oligomorphic.
Let $n \in {\mathbb N}$ be arbitrary. 
We have to verify that $\xi(G)$ has finitely
many orbits in its action on $N^n$. 
By (\ref{om-uniformity-1}), the equivalence relation
$V_n$ has finitely many equivalence classes $C_1,\dots,C_m$. Let $S,T \in N^n$ be
arbitrary, and let $s,t \in M^n$ be such that
$S_i = [s_i]_V$ and $T_i = [t_i]_V$ for all $i \in \underline n$. 
Suppose that $(s,t) \in C_j$ for some $j \in \underline m$. 
Then there exists a $g \in G$ such that $g(s_i,g t_i) \in V$ for all $i \in {\underline n}$. 
Hence, $\xi(g)([t_i]_V) = [gt_i]_V = [s_i]_V$,
and $S$ and $T$ lie in the same orbit
of the action of $\xi(G)$ on $N^n$. 
Therefore, $m$ bounds the number of orbits of $n$-tuples of $\xi(G)$. 

\vspace{.2cm}
{\bf Claim 3.} The map $\xi$ is continuous.
Let $x \in M$ be arbitrary, and $E \subseteq N$ be finite.
Define $W_{x,E} := \{g \in N^N  \; | \; g(e) = \xi(x)(e) \text{ for all $e \in E$} \}$, and recall from the proof of $(2) \Rightarrow (3)$ and $(3) \Rightarrow (4)$ that sets of the form $W_{x,E}$ form a 
uniformity base for the submonoid induced by $\xi(M)$ in $N^N$. 
Then 
\begin{align*}
\xi^{-1}(U_{x,E}) & = \{ y \in M \mid \xi(y)(e) = \xi(x)(e)\text{ for all $e \in E$} \} \\
& = \bigcap_{[z]_{V} \in E} \{ y \in M  \mid \xi(x)([z]_{V}) = \xi(y)([z]_{V})\}  \\
& = \bigcap_{[z]_{V} \in E} \{ y \in M  \mid [xz]_{V} = [yz]_{V} \} \\
& = \bigcap_{[z]_{V} \in E} \{ y \in M  \mid yz \in [xz]_{V} \} \, .
\end{align*}
We are going to prove that $S:=\{ y \in M  \; | \; yz \in [xz]_{V} \}$ is open.  
As $\bf M$ is a topological monoid, $\rho_z \colon M \to M, \, y \mapsto yz$ is continuous for every $z \in M$. Since $V$ is a member of a compatible uniformity on $M$, each element of $M/V$ is open in $M$, and in particular, $[xz]_{V}$ is open. Hence, $\rho_z^{-1}([xz]_{V}) = \{ y \in M \mid \rho_z(y) \in [xz]_{V} \} = S$ is open, too. Therefore, $\xi^{-1}(U_{x,E})$ is a finite intersection of open sets and therefore open. This shows that $\xi$ is continuous at $x$. 

\vspace{.2cm}
{\bf Claim 4.} The map $\xi^{-1}$ is uniformly continuous. According to (\ref{om-uniformity-2}), we have to show that for every finite subset $F \subseteq C$ the relation
$\{(\xi(x),\xi(y))  \mid (x,y) \in V_{F} \}$ 
is  an entourage in $\xi(M)$ where $V_{F} := \{ (x,y) \in M^2 \mid \forall g \in F \colon \, (xg,yg) \in V \}$. 
Define \begin{displaymath}
	U_F := \{(f,g) \in N^N \times N^N  \mid \forall z \in F \colon \, f([z]_{V}) = g([z]_{V}) \} .
\end{displaymath} Now, 
\begin{align*}
U_{F} \cap \xi(M)^2 = & \; \big \{(\xi(x),\xi(y)) \mid x,y \in M \text{ and } \forall g \in F \colon \, \xi(x)([g]_{V}) = \xi(y)([g]_{V}) \big \} \\
= & \; \big \{(\xi(x),\xi(y)) \mid x,y \in M \text{ and } \forall g \in F \colon \, [xg]_{V} = [yg]_{V} \big \}  \\ 
= & \; \big \{(\xi(x),\xi(y)) \mid (x,y) \in V_{F} \big \} \, . 
\end{align*}
The former set is an entourage in $\xi(M)$, and this proves the claim. 

\vspace{.2cm}
{\bf Claim 5.} The image of $M$ under $\xi$ is closed in $N^N$. This is analogous to the proof of the corresponding claim in Theorem~\ref{thm:monoids}, using uniform continuity of $\xi^{-1}$,
the fact that $\bf C$ has a compatible complete uniformity, and continuity of $\xi$. 
\end{proof}

\section{Oligomorphic Permutation Groups}
\label{sect:ogroups}
A topological group $\bf G$ is \emph{Roelcke precompact} if for every open 
subgroup $\bf U$ of $\bf G$ there
exists a finite set $F \subseteq G$
such that ${\bf U} F {\bf U} = G$. 
Todor Tsankov~\cite{Tsankov} showed that a closed
subgroup of $\Sym({\mathbb N})$ is Roelcke precompact 
if and only if every continuous action of $\bf G$ with finitely many orbits on a countable, discrete set $X$ 
is oligomorphic. 
This can be turned into a characterisation
of those subgroups of $\Sym({\mathbb N})$ that are topologically isomorphic to a closed oligomorphic permutation group on a countable set. 

\label{thm:ogroups}
\begin{theorem}\label{cor:ogroups}
Let $\bf G$ be a topological group. Then the following
are equivalent.
\begin{enumerate}
\item $\bf G$ is topologically isomorphic to the automorphism group of a countable $\omega$-categorical structure;
\item $\bf G$ is topologically isomorphic to a
closed oligomorphic subgroup of $\Sym({\mathbb N})$. 
\item $\bf G$ is Polish, Roelcke precompact, and has a compatible $G$-finitely generated left invariant ultrametric. 
\item $\bf G$ is Polish, Roelcke precompact, and has
an open 
subgroup $V$ of countably infinite index such that for all open subgroups $\bf U$ of $\bf G$ there are $g_1,\dots,g_n \in G$ such that $\bigcap_{i \in \underline n} g_i V g_i^{-1} \subseteq U$. 
\item $\bf G$ is topologically isomorphic to a
closed transitive oligomorphic subgroup of $\Sym({\mathbb N})$. 
\end{enumerate}
\end{theorem}
\begin{proof}
The equivalence of (1) and (2) 
is again well-known, and similar to 
the equivalence of (1) and (2) in Theorem~\ref{thm:omonoids}. 
For the implication from (2) to (3), 
we already 
have Roelcke precompactness of $\bf G$
from~\cite{Tsankov} as mentioned above. 
The group $\bf G$ is Polish as it is topologically isomorphic to a closed subgroup of
the Polish group $\Sym({\mathbb N})$. 
We now consider the monoid $\bf M$ induced 
in ${\mathbb N}^{\mathbb N}$ by the closure of the 
image of $G$ in ${\mathbb N}^{\mathbb N}$ under the topological isomorphism. By the implication
from (2) to (3) in Theorem~\ref{thm:omonoids} there 
exists a compatible $G$-finitely generated 
complete left non-expansive ultrametric.
Then the restriction of $d$ to $G$
is no longer guaranteed to be complete (since the image of $G$ may not be closed), but is still a compatible $G$-finitely generated left invariant ultrametric. 

\vspace{.2cm}
To show the implication from (3) to (4), 
let $d$ be a metric as in (3). 
Since $d$ is a left non-expansive metric, 
$V:=\{x \in G \mid d(x,1) < 1/2\}$ induces 
an open subgroup in $\bf G$. 
Let $\bf U$ be any open subgroup
of $\bf G$. Since $d$ is compatible and $\bf U$ is open,
there exists an $\varepsilon > 0$ such that
$\{x \in G \mid d(x,1) < \varepsilon\} \subseteq U$. 
Since $d$ is $G$-finitely generated 
there exists a finite subset $F$ of $G$
such that 
$\forall x,y \in G : \max_{g \in F} d(xg,yg) < 1/2 \Rightarrow d(x,y) < \varepsilon$. 
Let $x \in \bigcap_{g \in F} g V g^{-1}$.
Then for all $g \in F$ we have
$g^{-1}xg \in V$  and thus $d(g^{-1}xg,1) < 1/2$. 
Since $g$ is non-expansive we have $d(xg,g) < 1/2$, and the above then implies that
$d(x,1) < \varepsilon$, and $x \in U$. 
Hence, $\bigcap_{g \in F} g V g^{-1}  \subseteq U$, 
as desired. 


\vspace{.2cm}
The implication from (4) to (5) 
can be found in the proof of Theorem 7.3.1 in~\cite{Bodirsky-HDR}, and
the implication from (5) to (2) is trivial. 
\end{proof}

\section{Clones}\label{sect:clones}
Clones (in the literature often \emph{abstract clones}) relate to function clones in the same way as (abstract) groups relate to permutation groups: the elements of a clone correspond to the functions of a function clone, and the signature contains composition symbols to code how functions compose. Since a function clone 
contains functions of various arities, a clone will be formalised 
as a multi-sorted structure, with a sort for each arity. 

\begin{definition}
A \emph{clone} $\aC$ is a multi-sorted structure with sorts $\{C^{(i)} \; | \; i \in {\mathbb N}\}$ and the
signature $\{p_i^k \; | \; 1\leq i \leq k\}
\cup \{\comp^k_l \; | \; k,l \in {\mathbb N} \}$. The elements of the sort $C^{(k)}$ will be called the \emph{$k$-ary operations} of $\aC$. 
We denote a clone by 
$$\aC = (C^{(0)},C^{(1)},\dots;(p_i^k)_{1\leq i \leq k},(\comp^k_l)_{k,l \geq 1})$$ and require that
$p_i^k$ is a constant in $C^{(k)}$,
and that  
$\comp_l^k \colon C^{(k)} \times (C^{(l)})^k \to C^{(l)}$ is an
 operation of arity $k+1$.
Moreover, it holds that
\begin{align}
\comp_k^k(f,p_1^k,\dots,p_k^k) & = f \label{eq:pin} \\ 
\comp^k_l(p_i^k,f_1,\dots,f_k) & = f_i \label{eq:pout} \\
\comp_l^k(f,\comp_l^m(g_1,h_1,\dots,h_m),\dots,\comp_l^m(g_k,h_1,\dots,h_m)) 
& =  \nonumber  \\
\comp^m_l(\comp^k_m(f,g_1,\dots,g_k),h_1,\dots,h_m) & \, . \label{eq:comp}
\end{align}
\end{definition}
We write $C := \bigcup_{i \in \omega} C^{(i)}$,
and when $f \in C^{(k)}$ we call $\ar(f) := k$
the \emph{arity} of $f$.  

Every function clone $\sC$ gives rise to an abstract
clone $\aC$ in the obvious way: $p^k_i \in C^{(k)}$ denotes the projection $\pi^k_i \in \sC$,
and $\comp^k_l(f,g_1,\dots,g_k) \in C^{(l)}$ denotes
the composed function $$(x_1,\dots,x_l) \mapsto f(g_1(x_1,\dots,x_l),\dots,g_k(x_1,\dots,x_l)) \in \sC \, .$$
In the following, we will also use the term `abstract clone' in situations where we want to stress that we are 
working with a clone and \emph{not} with a function clone.

Let $\bf C$ and $\bf D$ be clones. 
A \emph{clone homomorphism} from $\bf C$
to $\bf D$ is a mapping $\xi$ from $C$ to $D$
such that 
\begin{itemize}
\item for all $i,k \in {\mathbb N}$, $i \leq k$,
the element $p_i^k$ in $C^{(k)}$ is mapped
to the element $p_i^k$ in $D^{(k)}$, and 
\item for all $f \in C^{(k)}$ and $g_1,\dots,g_k \in C^{(l)}$ 
$$\xi(\comp^k_l(f,g_1,\dots,g_k)) = \comp^k_l(\xi(f),\xi(g_1),\dots,\xi(g_k)) \, .$$
\end{itemize}
A \emph{(clone) isomorphism} is a bijective
clone homomorphism. 

\begin{definition}
A \emph{topological clone} is an abstract clone $\bf C$ together with a topology on $C$ such that the functions $(\comp^k_l)_{k,l \geq 1}$ of $\bf C$ are continuous.  
\end{definition}

A \emph{topological isomorphism} between topological
clones is an isomorphism between the respective 
abstract clones which is also a homeomorphism. 

Let $\bf C$ be a topological clone and $d$ a compatible metric on $C$. We say that $d$ 
is \emph{left non-expansive} if for all $f \in C^{(k)}$
and $g_1,h_1,\dots,g_k,h_k \in C^{(l)}$ we have that
$$d(\comp^k_l(f,g_1,\dots,g_k),\comp^k_l(f,h_1,\dots,h_k)) \leq \max(d(g_1,h_1),\dots,d(g_k,h_k)) \; .$$ 
We say that $d$ is \emph{projection right invariant} 
if for all $f,g \in C^{(k)}$, $k \in {\mathbb N}$, 
and pairwise distinct $i_1,\dots,i_l \in \underline k$ 
we have $$d(f,g) = d(\comp^l_k(f,p^k_{i_1},\dots,p^k_{i_l}),\comp^l_k(g,p^k_{i_1},\dots,p^k_{i_l}) \; .$$
There are several different ways to define 
compatible left non-expansive metrics for function clones. The following metric, 
however, is additionally projection right-invariant. 

\begin{definition}[Box metric]\label{def:box-m}
The \emph{box metric} on ${\mathbb N}^k \to {\mathbb N}$ 
is the metric
defined by $$d(f,g) = 2^{-\min \{ n \in \mN \, \mid \, \text{ there is } s \in \underline n^k \text{ such that } f(s) \neq g(s)\}} \,.$$ 
\end{definition}

A relation $R$ on $C$ is called 
\emph{left invariant} 
if for all $f \in C^{(k)}$ and $g_1,h_1,\dots,g_k,h_k \in C^{(l)}$, if $(g_1,h_1),\dots,(g_k,h_k) \in R$,  then $(\comp^k_l(f,g_1,\dots,g_k),\comp^k_l(f,h_1,\dots,h_k)) \in R$. 
A relation $R$ on $C$ is called
\emph{projection right invariant} if for all $k \in {\mathbb N}$ and pairwise distinct $i_1,\dots,i_l \in
 \underline k$
$$ (f,g) \in R \text{ implies }
(\comp^l_k(f,p^k_{i_1},\dots,p^k_{i_l}),\comp^l_k(g,p^k_{i_1},\dots,p^k_{i_l})) \in R \; .$$ 
As in the case of monoids, we might omit brackets
when composing functions in a clone; that is,
when $g \in C^{(l)}$ and $f \in (C^{(k)})^l$
we write $gf$ for $\comp^l_k(g,f)$.

\begin{thm}\label{thm:clones}
\label{satz:topo-clones}
Let $\bf C$ be a topological clone. Then the following are equivalent.
\begin{enumerate}
\item $\bf C$ is topologically isomorphic to 
the polymorphism clone of a countable structure.  
\item\label{c-start} $\bf C$ is topologically isomorphic to a closed subclone of ${\bf O}_{\mathbb N}$. 
\item\label{c-goal} $\bf C$ is separable and admits a compatible complete left non-expansive and projection right invariant ultrametric.
\item\label{uniformity} $\bf C$ is Hausdorff and has a complete compatible uniformity with a countable  base $\mathcal B$ such that each element of $\mathcal B$ is a left invariant and projection right invariant equivalence relation of countable index on $C$. 
\end{enumerate}
\end{thm}




\begin{proof}[Proof of Theorem~\ref{thm:clones}]
The equivalence between $(1)$ and $(2)$ is well-known; see, e.g.,~\cite{Szendrei}. 

$(2) \Rightarrow (3)$.
The proof that closed subclones of ${\bf O}_{\mathbb N}$ are separable is analogous to the
proof in the case of closed submonoids of ${\mathbb N}^{\mathbb N}$ given in Theorem~\ref{thm:monoids}. A compatible complete left non-expansive and projection right invariant ultrametric is given by the box metric (Definition~\ref{def:box-m}). 

$(3) \Rightarrow (4)$.
Let $d$ be a compatible left non-expansive and projection right invariant ultrametric of $\bf C$. Then $\mathcal B := \{B_n \; | \; n \in {\mathbb N} \}$ where $B_n := \{(x,y) \in C^2 \; | \; d(x,y) < 1/n\}$ is a countable uniformity base which we claim has the required properties. Since $d$ is complete, the uniformity is complete, too,
and it is clearly compatible with ${\bf C}$. 
The transitivity of the $B_n$ follows from
the assumption that $d$ is an ultrametric. 
Separability of ${\bf C}$ implies that the $B_n$
have countably many equivalence classes. Finally, 
$B_n$ is left invariant because $d$ is left non-expansive, and projection right invariant because $d$ is projection right invariant. 

$(4) \Rightarrow (2)$. 
Let $V_1,V_2,\dots$ be an enumeration of 
the equivalence relations in $\mathcal B$,
and let $(U_n)_{n \in \mathbb N}$ be the sequence of 
equivalence relations given by $U_n := \bigcap_{i \in \underline n} V_i$. Note that these relations are
left invariant and projection right invariant, too. 
Let $N_n := C/{U_n}$, and 
let $N$ be the countable set 
$\bigcup_{n \in \mathbb N} (N_n \times \{n\})$. 
Define $\xi \colon C \to {\mathscr O}_N$ as follows.
Let $f \in C^{(k)}$, $e_1,\dots,e_k \in N$, and for $i \in \underline{k}$ let $g_i \in C^{(l_i)}$ be such that 
$e_i = ([g_i]_{U_{n_i}},n_i)$  for $n_1,\dots,n_k \in {\mathbb N}$. 
Let $n := \min(n_1,\dots,n_k)$ and 
$l := \max(l_1,\dots,l_k)$. 
Define $$\xi(f)(e_1,\dots,e_k) := \big ([f(g_1(p^l_1,\dots,p^l_{l_1}),\dots,g_k(p^l_1,\dots,p^l_{l_k}))]_{U_n},n \big ) \, .$$

This is well-defined: 
let $h_1,\dots,h_k \in C$ 
be such that $(g_1,h_1) \in U_{n_1}, \dots, (g_k,h_k) \in U_{n_k}$. 
Then $(g_i(p^l_1,\dots,p^l_{l_i}),
h_i(p^l_1,\dots,p^l_{l_i})) \in U_{n_i}$  for all $i \in \underline k$ 
since $U_{n_i}$ is projection right invariant. 
Moreover, $(g_i(p^l_1,\dots,p^l_{l_i}),
h_i(p^l_1,\dots,p^l_{l_i})) \in U_n$
by the definition of $U_n$. 
It follows that $(f(g_1(p^l_1,\dots,p^l_{l_1}),\dots,g_k(p^l_1,\dots,p^l_{l_k})),f(h_1(p^l_1,\dots,p^l_{l_1}),\dots,h_k(p^l_1,\dots,p^l_{l_k})) \in {U_n}$ since $U_n$
is left invariant. 

\vspace{.2cm}
{\bf Claim 1.} The map $\xi$ is a homomorphism. 
Let $f \in C^{(k)}$ and $g_1,\dots,g_k \in C^{(l)}$.
Let $n_1,\dots,n_l \in {\mathbb N}$ 
and $h_1 \in C^{(m_1)},\dots,h_l \in C^{(m_l)}$. For $i \in \underline l$, 
define $e_i := ([h_i]_{U_{n_i}},n_i) \in N$. 
Let $m := \max(m_1,\dots,m_l)$ and let $n := \min(n_1,\dots,n_l)$. 
We write $h'_i$ for $h_i(p^m_1,\dots,p^m_{m_i})$. 
Then  
\begin{align*}
\xi(f(g_1,\dots,g_k))(e_1,\dots,e_l) = 
& \, ([f(g_1,\dots,g_n)(h'_1,\dots,h'_l)]_{U_n},n)  \\
= & \, ([f(g_1(h'_1,\dots,h'_l),\dots,
g_k(h'_1,\dots,
h'_l)))]_{U_n},n) \\
= & \, \xi(f)\big (([g_1(h'_1,\dots,h'_l)]_{U_n},n),
\dots,([g_k(h'_1,\dots,h'_l)]_{U_n},n) \big )\\
= & \, \xi(f)\big (\xi(g_1)(p_1,\dots,p_n),\dots,\xi(g_l)(e_1,\dots,e_n) \big )
\end{align*}
and
\begin{align*}
\xi(p^l_i)(e_1,\dots,e_l) 
= & \, ([p^l_i(h_1,\dots,h_l]_{U_n},n) \\
= & \, ([h_i]_{U_n},n) = e_i 
\end{align*}
and thus, $\xi(p^l_i)$ is mapped to the
 $i$-th $l$-ary projection in ${\mathscr O}_N$.

\vspace{.2cm}
{\bf Claim 2.} The map $\xi$ is injective. 
Let $f_1,f_2 \in C$ be distinct. As the topology on $C$ is Hausdorff 
there exists an $n \in {\mathbb N}$ such that
$(f_1,f_2) \notin U_n$.
Therefore, 
\begin{align*}
\xi(f_1)\big (([g_1]_{U_n},n),\dots,([g_k]_{U_n},n)\big ) = & \; [f_1(g_1,\dots,g_k)]_{U_n} \\
\neq & \; [f_2(g_1,\dots,g_k)]_{U_n} 
= \xi(f_2)\big (([g_1]_{U_n},n),\dots,([g_k]_{U_n},n) \big)
 \, .
\end{align*} 

\vspace{.2cm}
{\bf Claim 3.} The map $\xi$ is continuous. 
Let $f \in C^{(k)}$ be arbitrary, and $E \subseteq N$ be finite. 
Define $W_{f,E} := \{g \in \xi(C)  \; | \; g(e) = \xi(f)(e) \text{ for all } e \in E^k \}$, and observe that the sets of the form $W_{f,E}$ form a 
basis for the topology induced by $\xi(C)$ in ${\mathscr O}_N$. 
Then 
\begin{align*}
\xi^{-1}(W_{f,E}) & = \{ g \in C \; | \; \xi(g)(e) = \xi(f)(e)\text{ for all } e \in E^k \} \\
& = \bigcap_{e_1,\dots,e_k \in E} \big\{ g \in C  \; | \; \xi(g)(e_1,\dots,e_k) = \xi(f)(e_1,\dots,e_k)\big\}  
\end{align*}
Let $e_1,\dots,e_k \in E$ be arbitrary, 
and let $h_i \in C^{(l_i)}$ be such that 
$e_i = ([h_i]_{U_{n_i}},n_i)$. 
Let $n := \min(n_1,\dots,n_k)$ and 
$l := \max(l_1,\dots,l_k)$.  
Let $\rho_{h_1,\dots,h_k} \colon C^{(k)} \to C^{(l)}$ defined by $$\rho_{h_1,\dots,h_k}(f) := 
f(h_1(p^l_1,\dots,p^l_{l_1}),\dots,h_k(p^l_1,\dots,p^l_{l_1}))\; .$$ 
Now, 
\begin{align*}
S:= & \; \big \{ g \in C \; | \; \xi(g)(e_1,\dots,e_k) = \xi(f)(e_1,\dots,e_k) \big \} \\
= & \; \big \{ g \in C^{(k)}  \; | \; \rho_{h_1,\dots,h_k}(g) \in [\rho_{h_1,\dots,h_k}(f)]_{U_n} \big \} \\
= & \; \rho^{-1}_{h_1,\dots,h_k}([f]_{U_n}) \; .
\end{align*}
The set $[f]_{U_n}$ is open because $U_n$ is the member of a compatible uniformity. As composition in topological
clones is continuous, $\rho$ is continuous, 
and therefore $S$ is open, too. 
Hence, $\xi^{-1}(W_{f,E})$ is a finite intersection of open sets, and open. 
We have thus shown that for every open subset $V$ of $\xi(C)$
containing $\xi(f)$ there is an open subset of $C$
that contains $f$ and 
whose image is contained in $V$,
that is, $\xi$ is continuous at $f$. 

\vspace{.2cm}
{\bf Claim 4.} The map $\xi^{-1}$ is uniformly continuous. We have to show that for all $n \in \mathbb N$ the relation 
$\{(\xi(f),\xi(g))  \; | \; (f,g) \in U_n\}$ 
is  an entourage in $\xi(C)$. 
Define $W_n$ to be the set 
\begin{align*}
\big \{(f,g) \in {\mathscr O}_N \times {\mathscr O}_N  \; | \; & f( ([p_1^{\ar(f)}]_{U_n},n),\dots,([p^{\ar(f)}_{\ar(f)}]_{U_n},n)) \\
= & \; g(([p_1^{\ar(g)}]_{U_n},n),\dots,([p^{\ar(g)}_{\ar(g)}]_{U_n},n)) \big \}\, .
\end{align*}
Then 
\begin{align*}
W_n \cap \xi(C)^2 = & \; \big \{ (\xi(f),\xi(g)) \; | \; f,g \in C \text{ and } (f(p^{\ar(f)}_1,\dots,p^{\ar(f)}_{\ar(f)}), g(p^{\ar(f)}_1,\dots,p^{\ar(f)}_{\ar(f)})) \in U_n \big \} \\
= & \; \big \{(\xi(f),\xi(g)) \; | \; (f,g) \in U_n \big \} \, . 
\end{align*}
The former set is an entourage in $\xi(M)$, and this proves the claim. 

\vspace{.2cm}
{\bf Claim 5.} The image of $C$ under $\xi$ is closed in ${\mathscr O}_N$. 
Again analogous to the corresponding claim in the proof of Theorem~\ref{thm:monoids}. 
\end{proof}

\section{Oligomorphic Function Clones}
\label{sect:oligo-clones}
In this section we characterise those topological clones $\bf C$ that arise as polymorphism clones of countably infinite $\omega$-categorical structures. 
We cannot simply combine the requirements from the characterisation
of polymorphism clones (Theorem~\ref{thm:clones}) and the characterisation of oligomorphic transformation monoids for the unary part of $\bf C$ (Theorem~\ref{thm:omonoids}): in~\cite{Reconstruction}, a structure is described whose polymorphism clone does not have an action where
the units act transitively, but our proof
of Theorem~\ref{thm:omonoids} would
in this case produce a transitive action of the unary part of $\bf C$. 

Let $\bf C$ be a topological clone. Similarly to what we did in Section~\ref{sect:oligo-monoids}, we consider the group
$G := \{x \in C^{(1)} \mid \exists y \in C^{(1)} \colon \, xy = yx = p^1_1\}$ of its \emph{units}. Note that if $d$ is a left non-expansive compatible metric on $C$, then the left action $G \times C \to C, \, (g,x) \mapsto gx$ is an isometric action of $G$ on $(C,d)$. Furthermore, generalising Definition~\ref{def:g-finitely} from monoids to clones, we say that a metric $d$ on $C$ is \emph{$G$-finitely generated} if 
$$\forall \varepsilon > 0, k \in \mN \, \exists F \subseteq G \text{ finite } \forall x,y \in C^{(k)}: \sup_{g \in F^k} d(xg,yg) < 1/2 \Rightarrow d(x,y) < \varepsilon \, .$$

\begin{theorem}\label{thm:oclones}
Let ${\bf C}$ be a topological clone and let $G$ be the units of $\bf C$. Then the following are equivalent. 
\begin{enumerate}
\item ${\bf C}$ is topologically isomorphic to the polymorphism clone of an $\omega$-categorical structure. 
\item 
There is a topological isomorphism $\xi$ 
between $\bf C$ and a closed subclone
of ${\mathscr O}_{\mathbb N}$ such that 
$\xi(G)$ is oligomorphic. 
\item ${\bf C}$ is separable and admits a compatible $G$-finitely generated complete left non-expansive and projection right invariant ultrametric $d$ such that for every $k \in {\mathbb N}$ the left action of $G$ on $(C^{(k)},d)$ is approximately oligomorphic.
\item $\bf C$ is Hausdorff and there exists a left invariant and projection right invariant equivalence relation $V$ of countable index on $C$ such that 
\begin{enumerate}
\item the uniformity $\mathcal U$ generated by 
\begin{align*}
\big \{ \{(f,g) \in C^{(k)} \times C^{(k)} \mid k \in {\mathbb N} \text{ and } & \; \forall h \in F^k \colon \, (fh,gh) \in V\} 
\mid F \subseteq G \text{ finite} \big \}
\end{align*}
is compatible and complete;
\item for all $n,k \in {\mathbb N}$, the equivalence relation $V_{n,k}$ defined on $(C^{(k)})^{n}$  by
$$V_{n,k} := \big\{ (f,g) \mid  \exists h \in G \, \forall i \in \underline{n} \colon \,  (f_i, h g_i) \in V \big \}$$
has finite index.
\end{enumerate}
\end{enumerate}
\end{theorem}
\begin{proof}
The equivalence of (1) and (2) is a well-known consequence of the theorem of Engeler, Svenonius, and Ryll-Nardzewski and the equivalence of (1) and (2) in Theorem~\ref{thm:clones}.

\medskip
$(2) \Rightarrow (3)$. We have already seen in Theorem~\ref{thm:clones} that 
${\bf C}$ is separable and that it has a compatible complete left non-expansive and projection right invariant ultrametric $\tilde d$. Since $\xi(G)$ 
is an oligomorphic permutation group,
there exists some $n \in \mathbb N$
such that \begin{displaymath} 
	\mathbb N = \bigcup_{j \in \underline n} \{ \xi(g)(j) \mid g \in G\} .
\end{displaymath} Then $d$ defined by 
$d(x,y) := 2^n \tilde d(x,y)$ is also
a compatible left non-expansive and projection right invariant ultrametric on $M$. We prove that $d$ is furthermore 
$G$-finitely generated. 

Let $\varepsilon > 0$ and $k \in {\mathbb N}$. 
Pick $m \in {\mathbb N}$
such that $2^{n-m} < \varepsilon$. 
For each $i \in \underline m$, choose some $g_i \in G$ such that $i \in \xi(g)({\underline n})$. Evidently,
$F := \{g_i \; | \; i \in {\underline m}\}$ is finite. 
Now, let $x,y \in C^{(k)}$. If $\max_{g \in F^k} d(xg,yg) < 1/2$ then $\xi(xg)|_{\underline n} = \xi(yg)|_{\underline n}$ for each $g \in F^k$. 
Thus, $\xi(x)|_{{\underline m}^k} = \xi(y)|_{{\underline m}^k}$, which implies that 
$\tilde d(x,y) \leq 2^{-m}$, and hence
$d(x,y) \leq 2^{n-m} < \varepsilon$. 

To complete the proof of (3), 
arbitrarily choose $\varepsilon > 0$ and $k \in \mN$. 
Let $m \in \mN$ be such that $2^{n-m} < \varepsilon$. 
Since $\xi(G)$ is oligomorphic, it has finitely many orbits in its componentwise action 
on ${\mathbb N}^{nm^k}$. In the following, we identify the elements of $C$ with their image under $\xi$ when applying them to elements from $\mathbb N$. 
For $f \in (C^{(k)})^n$, we write $f[{\underline m^k}]$
for the $nm^k$-tuple 
\begin{align*}
\big (& f_1(1,\dots,1),f_1(1,\dots,1,2),\dots,f_1(m,\dots,m), \nonumber \\
& \dots \\ 
& f_n(1,\dots,1),f_n(1,\dots,1,2),\dots,f_n(m,\dots,m) \big ) \, . \nonumber
\end{align*}
Let $Q$ be the set of all orbits of $nm^k$-tuples
of the form $f[{\underline m}^k]$ for some 
$f \in (C^{(k)})^n$. For each $P \in Q$, choose an $f_P \in (C^{(k)})^n$ such that $f_P[{\underline m}^k] \in P$. Clearly, the 
set $F := \{f_P \; | \; P \in Q\}$ is a finite subset of 
$(C^{(k)})^n$. Let $f \in (C^{(k})^n$ be arbitrary.
Then there exists some $P \in Q$ such that $f[{\underline m}^k] \in P$. Hence, there exists $g \in G$ such that $f[{\underline m}^k] = g(f_P[{\underline m}^k])$. Consequently, for every $j \in \underline n$ we have $\tilde d(f_j,f_{P,j}) \leq 2^{-m}$ and hence 
$d(f_j,f_{P,j}) \leq 2^{n-m} < \varepsilon$. 

\medskip 
$(3) \Rightarrow (4)$. 
Let $d$ be as in (3). Since $d$ is an ultrametric, $V := \{(x,y) \in C^2 \mid d(x,y) < 1/2\}$ is an equivalence relation. As $d$ is left non-expansive, 
$V$ is left-invariant. 
As $d$ is projection right invariant, 
$V$ is projection right invariant. 
Moreover, $C/V$ is countable because $C$ is separable. 
We now show (4)(a). 
By definition, the uniformity $\mathcal U$ is generated by 
\begin{align*}
\big \{ \{(x,y) \in C^{(k)} \times C^{(k)} \mid k \in {\mathbb N} \text{ and } \forall g \in F :
d(xg,yg) < 1/2 \} \mid F \subseteq G \text{ finite} \big \}
\; .
\end{align*} 
Since $d$ is $G$-finitely generated we find for all $n,k \in \mathbb N$ a finite set $F \subseteq (C^{(k)})^n$ such that for all
$x,y \in (C^{(k)})^n$ we have that
$$\max_{g \in F} d(xg,yg) < 1/2 \Rightarrow d(x,y) < 1/n \; .$$
It follows that $\mathcal U$ is generated
by $$\big \{ \{(x,y) \in C^2 \mid d(x,y) < 1/n\} \mid n \in \mathbb N \big \} \; .$$
and hence
coincides with the uniformity generated by $d$. In particular, $\mathcal{U}$ is complete and compatible with the topology of $C$.

We finally show (4)(b). 
Let $n,k \in \mN$. By 
(3) there exists a finite set $F \subseteq (C^{(k)})^n$ such that for all $x \in (C^{(k}))^n$ there
exists a $y \in F$ such that
$d_{G,n}(x,y) < 1/n$. 
Hence, 
the index of $V_{n,k} = \big \{(x,y) \in C^{(k)} \times C^{(k)} \mid d_{G,n}(x,y) < 1/n \big \}$ 
is bounded by $|F|$.

\medskip
$(\ref{m-uniformity}) \Rightarrow (\ref{m-start})$.
By assumption, $N := C/V$ is countable. 
Let us define $\xi \colon C \to {\mathscr O}_N$
as follows. Let $f \in C^{(k)}$, 
$g_1 \in C^{(l_1)}, \dots, g_k \in C^{(l_k)}$, and $l := \max(l_1,\dots,l_k)$. 
Define
$$\xi(f)([g_1]_V,\dots,[g_k]_V) := [f(g_1(p_1^l,\dots,p^l_{l_1}),\dots,g_k(p_1^l,\dots,p_{l_k}^l))]_V \; .$$ 
We first show that the function $\xi$ is well defined. Let $h_1,\dots,h_k \in C$ be such that $(g_i,h_i) \in V$ for all $i \in \underline{k}$. 
Then 
$(g_i(p^l_1,\dots,p^l_{l_i}),
h_i(p^l_1,\dots,p^l_{l_i})) \in V$ since $V$ is projection right invariant. 
It follows that $(f(g_1(p^l_1,\dots,p^l_{l_1}),\dots,g_k(p^l_1,\dots,p^l_{l_k})),f(h_1(p^l_1,\dots,p^l_{l_1}),\dots,h_k(p^l_1,\dots,p^l_{l_k}))) \in V$ since $V$ is left invariant. 

\vspace{.2cm}
{\bf Claim 1.}
The map $\xi$ is a homomorphism.  
Let $f \in C^{(k)}$ and $g_1,\dots,g_k \in C^{(l)}$. 
Let $h_1,\dots,h_l \in C$ be arbitrary,
let $m_i$ be the arity of $h_i$ for $i \in \underline l$, and let $m := \max(m_1,\dots,m_l)$. 
We write $h'_i$ for $h_i(p^m_{1},\dots,p^m_{m_i})$. 
Then 
\begin{align*}
\xi(f(g_1,\dots,g_k))([h_1]_V,\dots,[h_l]_V) 
= & \; [f(g_1,\dots,g_k)(h'_1,\dots,h'_l)]_V \\
= & \; [f(g_1(h'_1,\dots,h'_l),\dots,g_k(h_1',\dots,h_l'))]_V \\
= & \; \xi(f)([g_1(h_1',\dots,h_l')]_V,\dots,[g_k(h_1',\dots,h'_l)]_V) \\
= & \; \xi(f)(\xi(g_1)([h_1]_V,\dots,[h_l]_V),\dots,\xi(g_k)([h_1]_V,\dots,[h_l]_V))
\end{align*} 
and $\xi(p^l_i)([h_1]_V,\dots,[h_l]_V) = [p^l_i(h_1,\dots,h_l)]_V = [h_i]_V$.

\vspace{.2cm}
{\bf Claim 2.} The map $\xi$ is injective. 
Let $f,g \in C$ be distinct. If $f$ and $g$
have distinct arities, then clearly $\xi(f) \neq \xi(g)$, so suppose that $\ar(f)=\ar(g)=k \in \mathbb N$. 
Due to (4b) and $C$ being Hausdorff, there are $h_1,\dots,h_k \in G$ such that $(f(h_1,\dots,h_k),g(h_1,\dots,h_k)) \notin V$.
Accordingly, $$\xi(f)([h_1]_V,\dots,[h_k]_V) = [f(h_1,\dots,h_k)]_V \ne [g(h_1,\dots,h_k)]_V = \xi(g)([h_1]_V,\dots,[h_k]_V)$$
and thus $\xi(f) \neq \xi(g)$. 

\vspace{.2cm}
{\bf Claim 3.} $\xi(G)$ is oligomorphic. 
The proof is as in the proof of Theorem~\ref{thm:omonoids}. 

\vspace{.2cm}
{\bf Claim 4.} The map $\xi$ is continuous.
Let $f \in C^{(k)}$ be arbitrary, and $F \subseteq N$ be finite.
Define $W_{f,F} := \{g \in {\mathscr O}_N  \; | \; g(e) = \xi(f)(e) \text{ for all $e \in F^k$} \}$, and recall that the sets of the form $W_{f,F}$ form a 
uniformity base for the subclone induced by $\xi(C)$ in ${\mathscr O}_N$. 
Then 
\begin{align*}
\xi^{-1}(U_{f,F}) & = \{ y \in C \mid \xi(y)(e) = \xi(f)(e) \text{ for all $e \in F^k$} \} \\
& = \bigcap_{[h_1]_V,\dots,[h_k]_V \in F} \big \{ g \in C  \mid [f(h_1,\dots,h_k)]_V = [g(h_1,\dots,h_k)]_V \big \} \\
& = \bigcap_{[h_1]_V,\dots,[h_k]_V \in F} \big \{ g \in C  \mid g(h_1,\dots,h_k) \in [f(h_1,\dots,h_k)]_V \big \} \, .
\end{align*}
We prove that $S :=\{ g \in C  \mid g(h_1,\dots,h_k) \in [f(h_1,\dots,h_k)]_V \}$ is open.  
As $\bf C$ is a topological clone, $\rho_{h_1,\dots,h_k} \colon C^{(k)} \to C^{(l)}, \, g \mapsto g(h_1,\dots,h_k)$ is continuous for all $(h_1,\dots,h_k) \in (C^{(l)})^{k}$. 
Since $V$ is both an equivalence relation and a member of a compatible uniformity on $C$, it follows that each element of $C/V$ is open in $M$, and, in particular, $[f(h_1,\dots,h_k)]_V$ is open. 

Hence, $\rho_{h_1,\dots,h_k}^{-1}([f(h_1,\dots,h_k)]_V) = \{ g \in C \mid \rho_{h_1,\dots,h_k}(g) \in [f(h_1,\dots,h_k)]_V \} = S$ is open, too. Therefore, $\xi^{-1}(U_{f,F})$ is a finite intersection of open sets and therefore open. This shows that $\xi$ is continuous at $f$. 

\vspace{.2cm}
{\bf Claim 5.} The map $\xi^{-1}$ is uniformly continuous. According to (4b), we have to show that for every finite subset $F \subseteq G$ the relation
$\{(\xi(f),\xi(g))  \mid (f,g) \in V_F \}$ 
is  an entourage in $\xi(C)$ where $V_F := \{ (f,g) \in C^2 \mid \forall h \in F^k \colon \, (fh,gh) \in V \}$. 
Define \begin{displaymath}
	U_F := \bigcap_{h_1,\dots,h_k \in F} \{(f,g) \in {\mathscr O}_N \times {\mathscr O}_N  \mid f([h_1]_V,\dots,[h_k]_V) = g([h_1]_V,\dots,[h_k]_V) \big \} .
\end{displaymath} Now, 
\begin{align*}
& U_F \cap \xi(C)^2 \\
= & \; \big \{(\xi(f),\xi(g)) \mid f,g \in C \text{ and } \forall h_1,\dots,h_k \in F \colon \, \xi(f)([h_1]_V,\dots,[h_k]_V) = \xi(g)([h_1]_V,\dots,[h_k]_V) \big \} \\
= & \; \big \{(\xi(f),\xi(g)) \mid f,g \in C \text{ and } \forall h_1,\dots,h_k \in F \colon \, [f(h_1,\dots,h_k)]_V = [g(h_1,\dots,h_k)]_V \big \}  \\ 
= & \; \big \{(\xi(f),\xi(g)) \mid (f,g) \in V_F \big \} \, . 
\end{align*}
The former set is an entourage in $\xi(C)$, and this proves the claim. 

\vspace{.2cm}
{\bf Claim 6.} The image of $C$ under $\xi$ is closed in ${\mathscr O}_N$. 
This is analogous to the corresponding claim in the proof of Theorem~\ref{thm:monoids}. 
\end{proof}


\section{Open Problems}
There is an important property of oligomorphic permutation groups, transformation monoids, and function clones that only depends on the topological group, topological monoid, and topological clone, respectively:
we say that an oligomorphic permutation group $\mathscr G$ (oligomorphic transformation monoid $\mathscr M$, oligomorphic function clone $\mathscr C$) is \emph{finitely related}
if there exists a countable structure $\Gamma$ with finite relational signature such that ${\mathscr G} = \Aut(\Gamma)$ (${\mathscr M} = \End(\Gamma)$, ${\mathscr C} = \Pol(\Gamma)$). 
It has been shown in~\cite{Topo-Birk} that all oligomorphic function clones that are
topologically isomorphic to a finitely related oligomorphic function clone must also be finitely related. 
The same argument also shows the corresponding statements for oligomorphic transformation monoids
and oligomorphic groups.  
Hence, it is natural to ask for a direct \emph{topological}
characterisation of finite relatedness in these three
settings, similarly to the characterisations obtained 
in this paper. 

\section*{Acknowledgements}
We would like to thank the referee for very carefully reading
our journal submission. 

\bibliographystyle{alpha}
\bibliography{local.bib}

\end{document}